\documentclass[11pt]{article}
\usepackage[margin=1in]{geometry}
\usepackage{amsmath,amsfonts,amssymb}
\usepackage[font=normalsize,singlelinecheck=off,justification=centerlast]{subcaption}
\usepackage{graphicx}
\usepackage{gensymb}

\usepackage{hyperref}
\usepackage{cleveref}
\usepackage{float}
\usepackage{pifont}
\usepackage{authblk}
\usepackage[style=numeric,backend=bibtex]{biblatex}
\usepackage{abstract}
\usepackage{xcolor}

\newcommand{\revans}[1]{\textcolor{black}{#1}}
\newcommand{\xx}{\mathbf{x}}

\begin{document}

\title{Numerical simulation of transient heat conduction with moving heat source using Physics Informed Neural Networks}

\author[1]{Anirudh Kalyan}
\author[1]{Sundararajan Natarajan\thanks{Corresponding author: snatarajan@iitm.ac.in}}

\affil[1]{
Indian Institute of Technology Madras\\
}

\date{June 5, 2025}

\maketitle
\vspace{-12px}

\begin{abstract}
In this paper, the physics informed neural networks (PINNs) is employed for the numerical simulation of heat transfer involving a moving source. To reduce the computational effort, a new training method is proposed that uses a continuous time-stepping through transfer learning. Within this, the time interval is divided into smaller intervals and a single network is initialized. On this single network each time interval is trained with the initial condition for $(n+1)^{\rm th}$ as the solution obtained at $n^{\rm th}$ time increment. Thus, this framework enables the computation of large temporal intervals without increasing the complexity of the network itself. The proposed framework is used to estimate the temperature distribution in a homogeneous medium with a moving heat source. The results from the proposed framework is compared with traditional finite element method and a good agreement is seen.
\end{abstract}

\vspace{10px}

\section{Introduction}\label{sec0}

\begin{figure}[htpb]
\subfloat[With time stepping (our method)]{\includegraphics[scale=0.30]{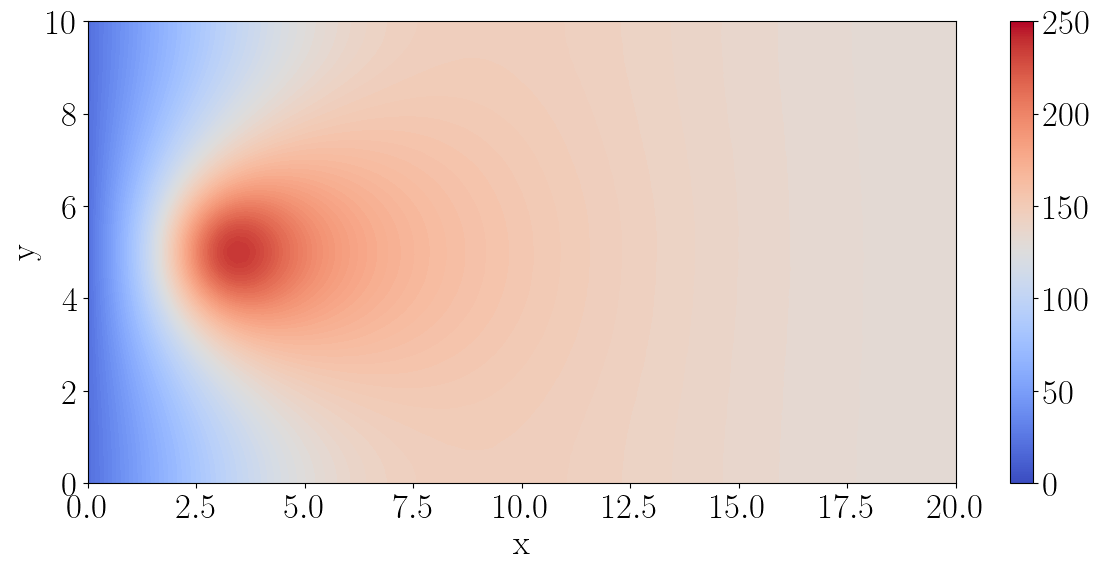}}
\subfloat[Without time stepping]{\includegraphics[scale=0.30]{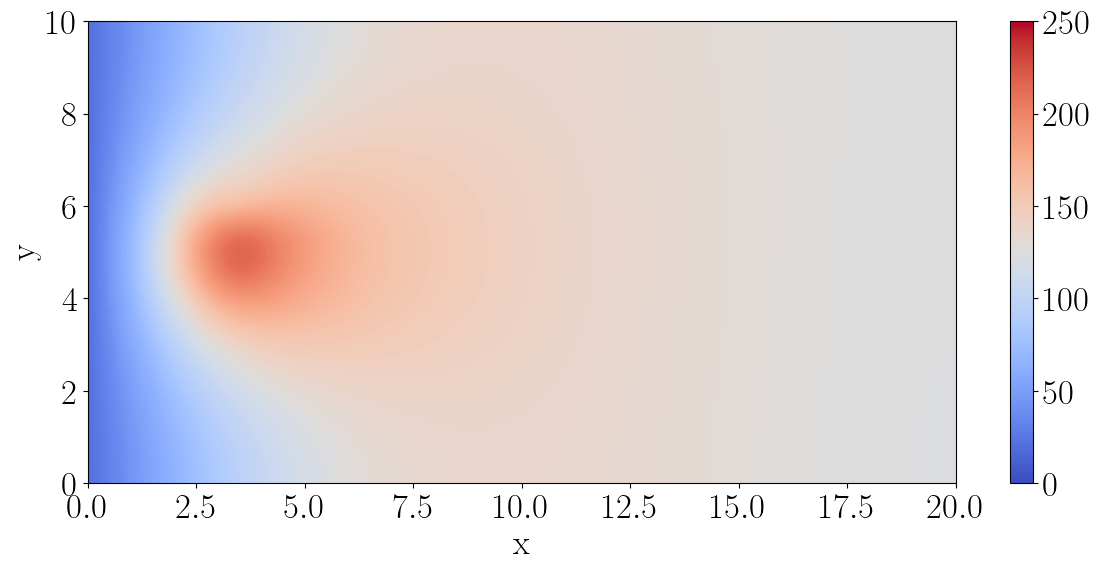}}
\caption{Temperature distribution output of PINN at $t=8s$ (see \cref{fig:geom_des}) with and without using time-stepping}
\label{fig:lineplotcompa_velocity}
\end{figure}

For the manufacturing processes such as grinding, welding, additive manufacturing , laser and plasma cutting , the physical phenomenon are governed by transient heat conduction. Amongst available numerical approaches, the finite element method (FEM) is the most commonly adopted approach. The FEM relies on dividing the domain into nonoverlapping subdomains, called the mesh connected at points, referred to as nodes. However, despite its wider application, it has the following drawbacks: (a) the quality of the solution is highly influenced by the quality of the underlying discretization and (b) because it uses linear piecewise approximations for representing the unknown fields, it requires very fine discretization to capture steep thermal gradients, especially near the heat source~\cite{Mundra01021996}. The latter can be addressed by adopting mesh refinement techniques~\cite{kuangatluri1985,kim2012} or using higher-order approximations. However, implementing adaptive mesh refinement technique along the heat source is challenging, especially in three dimensional problems.

Rosenthal first applied the heat diffusion equation for a moving source\cite{rosenthal1941mathematical}, leading to the Rosenthal formulation of the heat diffusion equation. Kunag and Atluri~\cite{kuangatluri1985} employed a moving mesh technique to estimate the temperature distribution due to a moving heat source. Kim~\cite{kim2012} used adapted mesh generation technique that ensures finer discretization in the path of the heat source to solve a heat conduction with a moving heat source. Storti et al.,~\cite{stortialbanesi2022} proposed a Chimera-FEM for transient heat conduction problem with moving heat source. Within this framework, sharp gradients are captured by fine mesh without requiring a globally fine fixed mesh or adaptive mesh refinement. Res\'endiz-Flores and Saucedo-Zendejo~\cite{resendizfloreszendejo2015} proposed a finite point set method to simulate temperature distribution in welding process. The salient feature of this is that the method works directly with the strong form of the governing partial differential equation (PDEs). Hostos et al.,~\cite{hostosstorti2022} proposed an improved element-free Galerkin method for a moving heat source problem in arc welding process. It uses two overlapping distribution of nodes for solving the thermal problem and the transfer of temperature distribution and heat flux is achieved by properly defined immersed boundaries. Mohammadi et al.,\cite{mohammadihematiyan2021} proposed a boundary-type meshfree method for 2D transient heat conduction with moving heat sources. It is based on method of fundamental solutions and the governing equations and boundary conditions are satisfied exactly. The introduction of mesh-free techniques~\cite{hostosstorti2022} and isogeometric analysis~\cite{zangliu2021,gongqin2022} relaxes the constraint imposed by the FEM on the mesh topology. Gong et al.,~\cite{gongqin2022} proposed an isogeometric boundary element analysis to solve heat conduction problems in multiscale electronic structures. The distinctive advantage of the approach is that only the boundary of the domain is discretized and the method benefits from computational efficiency, thanks to the isogeometric analysis. The problem of accurately modelling steep thermal gradients is partially addressed by augmenting the finite element basis with anastz within the framework of the extended finite element method~\cite{muntshulshoff2003,haraduarte2009,bordasmenk2024}. For a comprehensive overview of different computational strategies for manufacturing based on laser bed fusion, interested readers are referred to~\cite{sarkarkapil2024}. However, despite these advances, there is still growing interest in improving the numerical approaches to model the moving heat source problem.

Recently, an alternate approach that utilizes Artificial Neural Networks (ANNs) as local basis for representing the unknown fields has gained growing interest amongst researchers. Based on the universal approximation theorem~\cite{HORNIK1989359} that states that any continuous function can be represented by a single layer network with an arbitrary number of neurons. The salient feature of the ANNs is its superior capability to represent complex non-linear data~\cite{ling_kurzawski_templeton_2016, PARISH2016758}. The idea of using neural networks (NN) to obtain an approximate solution of governing differential equations, ordinary and partial differential equations (PDEs) by modifying the loss function of ANNs was demonstrated in~\cite{712178,870037}. This was further propelled by Karniadakis and co-workers~\cite{RAISSI2019686} in a series of publications leading to what is now popularly called as the `Physics Informed Neural Networks (PINNs). The fundamental idea behind the PINN is to approximate the solution to the PDE by Deep Neural Network and convert it to solving a nonlinear optimization problem. It employs automatic gradient computations~\cite{10.5555/3122009.3242010} and typically are implemented using Machine-learning libraries, e.g., TensorFlow \cite{45381}. \revans{This simplicity has led to an exponential growth of the application of PINNs to approximate the solution of PDEs~\cite{fuxu2024,gaofu2024,tangfu2022,cai2021physics,cai2021physicsHT,haghighat2020deep,dwivedi2019distributed}. PINNs have been used to solve a wide variety of differential equations including fluid mechanics\cite{cai2021physics}, heat transfer\cite{cai2021physicsHT}, solid mechanics\cite{haghighat2020deep}, solutions to PDEs on surfaces~\cite{tangfu2022} and thermo-hydro coupled problems in concrete~\cite{gaofu2024}}. Dwivedi et al.,~\cite{dwivedi2019distributed, dwivedi2019physics} introduced domain-decomposed version of PINN, called distributed PINN (DPINN) where an additional ``interface-loss" term is added in order to account for continuity and differentiability requirements at interfaces. A seamless integration between the traditional methods and the PINNs was proposed by Ramabathiran anad Ramachandran~\cite{RAMABATHIRAN2021110600} by using DNNs as kernels. Goswami et al.,~\cite{GOSWAMI2020102447} developed variational-energy-based PINN to address the problem of brittle fracture. Energy-based approaches were also proposed for elastic problems~\cite{SAMANIEGO2020112790, ZHUANG2021104225, NGUYENTHANH2021114096}. Herrmann and Kollmannsberger~\cite{hermannkollmannsberger2024} provided an overview of the application of deep learning techniques to computational mechanics. A comprehensive overview of application of artificial intelligence for solving PDEs is discussed by Wang et al.,~\cite{wangbai2024}. The paper discusses different approaches such as Physics-Informed Neural Networks (PINNs), Deep Energy Methods (DEM), Operator Learning, and Physics-Informed Neural Operator (PINO) for solving the PDEs. \revans{Recently, in a series of papers, Liu and co-workers~\cite{qiuwang2025,qiuwang2024,zhangwang2023} have adopted the framework of PINNs for multi-physics problems for homogeneous and heterogeneous medium, including functionally graded materials.}

To the best of author's knowledge, the PINNs have not been applied to transient heat conduction problem with a moving heat source. Given the wider applicability and simplicity, we propose to use PINNs to achieve a truly meshless approach to solve the transient heat conduction with a moving heat source. The novelty of the current work includes: (a) application of the PINNs to moving heat source problem and (b) a new training method that uses continuous time-stepping through transfer learning. Within this new training method, the solution at $n^{\rm th}$ time step is used as initial condition for the $(n+1)^{\rm th}$ time step. Thus, enabling the computation of large temporal intervals without increasing the complexity of the network itself.

The rest of the manuscript is organized as follows: \Cref{sec:pinnoverview} presents an overview of physics informed neural networks. The governing equations for a heat transfer with moving heat source is presented in \Cref{sec:theory}.  The section also discusses the application of PINNs to the transient heat conduction problem. In \Cref{sec:numexa}, we demonstrate the accuracy of the PINNs for transient heat conduction with a moving heat source, the results are compared with traditional finite element method. We conclude the paper with a brief summary and outlook in \Cref{sec:concl}.

\section{Overview of the physics informed neural network}
\label{sec:pinnoverview}
Physics Informed Neural Networks (PINNs) represent a new paradigm for solving differential equations that leverage recent developments in parallel processing and deep learning techniques. Central to the operation of PINNs is the universal approximation theorem~\cite{HORNIK1989359}, which states that any continuous function can be approximated with a network of neurons with weight adjusted to minimize the approximation error.

A typical neural network (See Figure 1) itself is composed of neurons, which contain 2 variables within themselves, weights $W$ and biases $b$. Each neuron represents a calculation $Wx+b$, where $x$ is the input to the neuron. The output of each neuron is passed through a non-linear activation function $\sigma$. Hence, the output of a single neuron is represented as $\sigma(Wx+b)$
\begin{figure}[hptb]
\centering 
\includegraphics[scale=0.35]{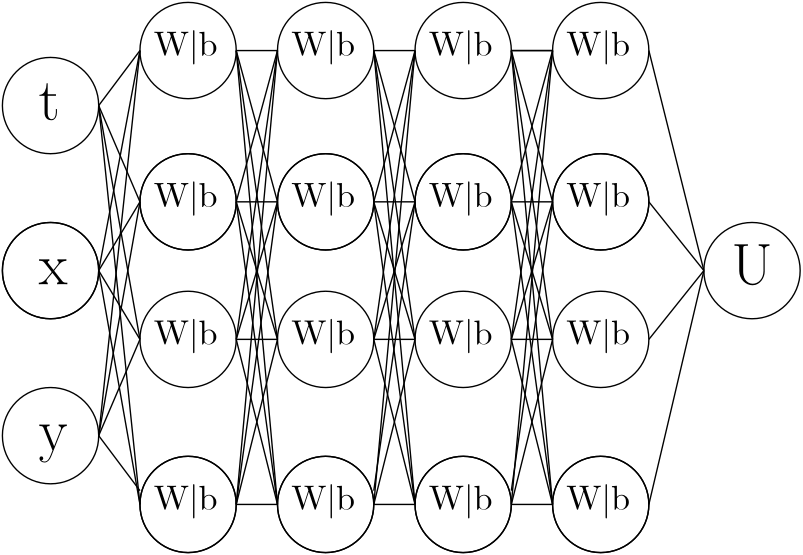}
\caption{Schematic representation of four layer neural network: $W,b$ represent the weights and biases of each neuron, $t,x,y$ are the input variables and $U$ is the associated output variable.}
\label{fig:nnet}
\end{figure}
Multiple neurons are connected to form a deep neural network. To train a network, the concept of back-propagation is utilized. Back-propagation utilizes the concept of chain rule to propagate the derivatives back to the inputs. Modern autograd mechanisms built into libraries like JaX\cite{bradbury2021jax}, and PyTorch\cite{paszke2017automatic} automate the chain rule process by using tensor graphs and parallelizing the procedure using GPUs, significantly reducing training times. 

Now, let us consider the following differential equation: 
\begin{equation}
    \dfrac{\partial u}{\partial t} + \mathcal{N}[u] = 0,
    \label{eqn:genericpde}
\end{equation}
where, $\mathcal{N}[u]$ denotes a nonlinear spatial differential operator applied to the function $u$. Let $\theta$ represent parameters of a network with weights $W$ and biases $b$. To find an approximate solution $u_\theta$ that satisfies \Cref{eqn:genericpde}, the PINN utilizes a loss function, which is a weighted sum of three terms, namely: (a) loss of the governing differential equation, $\mathcal{L}_{r}$; (b) loss corresponding to the boundary conditions (Neumann or Dirichlet), $\mathcal{L}_{bc}$ and (c) loss corresponding to the initial condition, $\mathcal{L}_{ic}$, given by:
\begin{equation}
\mathcal{L}(\theta) = \lambda_{ic} \mathcal{L}_{ic}(\theta) + \lambda_{bc} \mathcal{L}_{bc}(\theta) + \lambda_{r} \mathcal{L}_{r}(\theta).
\end{equation}
where $\lambda_i, i=ic, bc, r$ are the associated weights and
\begin{equation}
\begin{aligned}
&\mathcal{L}_r(\theta) = \frac{1}{N_r} \sum_{i=1}^{N_r} \left\| \frac{\partial u_{\theta}}{\partial t}(\xx,t) + \mathcal{N}[u_{\theta}](\xx,t) \right\|_2^2 \\
&\mathcal{L}_{ic} = \frac{1}{N_i} \sum_{i=1}^{N_i} \left[ u(\xx_i,0) - u_{0}(\xx_i) \right]^2 \\
&\mathcal{L}_{bc} = \frac{1}{N_b} \sum_{i=1}^{N_b} \left[ u(\xx_{b_i},t_i) - u_{b}(\xx_{b_i},t_i) \right]^2 \quad {\rm For~Dirichlet~boundary~conditions} \\
&\mathcal{L}_{bc} = \frac{1}{N_b} \sum_{i=1}^{N_b} \left[ \frac{\partial u}{\partial n}(\xx_{b_i},t_i) - \frac{\partial u_b}{\partial n}(\xx_{b_i},t_i) \right]^2 \quad {\rm For~Neumann~boundary~conditions}
\end{aligned}
\label{eqn:lossfunction_def}
\end{equation}
where, $u_b$ is the Dirichlet boundary condition values and $\dfrac{\partial u_b}{\partial n}$ is the specified Neumann condition, respectively, at coordinates $(\xx_b)$, $n$ is the unit normal vector, $u_o$ is the data supplied to the network with respect to the initial conditions, $N_r$ is the total number of collocation points, which are the points at which the differential equation loss is enforced. $N_i$ is a the set of points where initial conditions are specified and $N_b$ represents the set of points where Dirichlet or Neumann conditions are specified. Within the framework of PINNs, the computational domain is initialized with a distribution of collocation points and the loss term corresponding to the differential equation $\mathcal{L}_r(\theta)$ is computed, i.e., the strong form (c.f. \Cref{eqn:genericpde}) is directly used. The error in the loss function is minimized in a least-square sense. Other implementations of PINNs incorporate the weak form of the differential equation\cite{kharazmi2019variational}\cite{kharazmi2021hp}.

PINNs achieve machine-level accuracy in derivatives by taking advantage of automatic differentiation, which uses the chain rule to propagate the derivatives from the output layer up to the input layer. This enables derivatives of the outputs with respect to the inputs, multiple times without facing numerical instability as is the case with finite difference. Libraries like JaX enable forward mode automatic differentiation~\cite{margossian2019review}, which facilitates faster computation, reducing training times. The number of collocation points and the design of the network both play a crucial role in determining the accuracy and resolution of the solution, similar to how meshing and time steps play a role in FEM.

\section{Application of PINNs to transient heat conduction}
\label{sec:theory}

\subsection{Governing equations}
Consider a heterogeneous conducting body occupying $\Omega \subset \mathbb{R}^d,~d=2,3$ bounded by $\Gamma$, a $d-1$ surface. The boundary accommodates the following decomposition: $\Gamma = \Gamma_D \bigcup\Gamma_N$ and $\Gamma_D \bigcap \Gamma_N = \emptyset$, where $\Gamma_N$ and $\Gamma_D$ represent the Neumann and Dirichlet part of the boundary. The general governing differential equation with the boundary conditions for a non-linear transient temperature field is given by:
\begin{equation}\label{eq:def}
\begin{aligned}
\gamma(\xx) \frac{\partial u(\xx,t)}{\partial t}-\boldsymbol \nabla \cdot (k(\xx) \boldsymbol \nabla u(\xx,t)) =f(\xx,t) \quad {\rm on}  \quad \Omega \times (0,T]\\
u(\xx,t)=\hat{u} \quad {\rm in}  \quad \Gamma_D \times (0,T]\\
 k(u)\boldsymbol \nabla u \cdot \mathbf{n} =\mathbf{\hat{q}} \quad {\rm in}  \quad \Gamma_N \times (0,T]\\
u(\xx,t=0)=\hat{u_0} \quad {\rm in}  \quad \Gamma_D |_{t=0}\\
\end{aligned}
\end{equation}
\noindent
where $\boldsymbol \nabla=[\frac{\partial}{\partial x}\quad \frac{\partial}{\partial y}]^T$ is the differential operator, $u(\xx,t)$ is the temperature field to be determined, and $\mathbf{\hat{q}}$ is the heat flux at boundary $\Gamma_N$. $k(\xx)$ is the thermal conductivity, while,  $\gamma(\xx)=\rho(\xx)C(\xx)$, where $\rho(\xx)$ is the density and $C(\xx)$ is the specific heat constant. For the current numerical study, the heat source is modeled as a Gaussian:
\begin{equation}
    f(\xx,t) = Q_0 e^{-r^2 / r_0^2}
\end{equation}
where, $Q_0$ represents the maximum energy power, $r$ denotes the radial distance from the center of the heat source, and $r_0$ is the radius of the heat source.

\subsection{PINNs for transient heat conduction}
Unlike the traditional FEM, which converts the strong form of the differential equation given in \Cref{eq:def} into an equivalent weak form, the proposed framework utilizes the strong form directly as described in \Cref{sec:pinnoverview}. The loss function that corresponds to the governing differential equation is given by:
\begin{equation}
    \mathcal{L}_r(\theta) = \frac{1}{N_r} \sum_{i=1}^{N_r} \left\| \gamma(\xx) \frac{\partial u_{\theta}(\xx,t)}{\partial t}-\boldsymbol \nabla \cdot (k(\xx) \boldsymbol \nabla u_{\theta}(\xx,t)) - f(\xx,t) \right\|_2^2
\end{equation}
Depending on the initial conditions and the boundary conditions, corresponding loss functions (c.f. \Cref{eqn:lossfunction_def}) are used. With this, all the constraints pertaining to the network are specified.

Since the network maps $u(\xx,t)$, in space and time, scaling the network with time proves to be costly in terms of parameters required and training times, proving extrapolation with respect to time, computationally demanding. To mitigate this, a new sequential training regime is introduced. 
\begin{figure}[hptb]
\centering 
\includegraphics[scale=0.35]{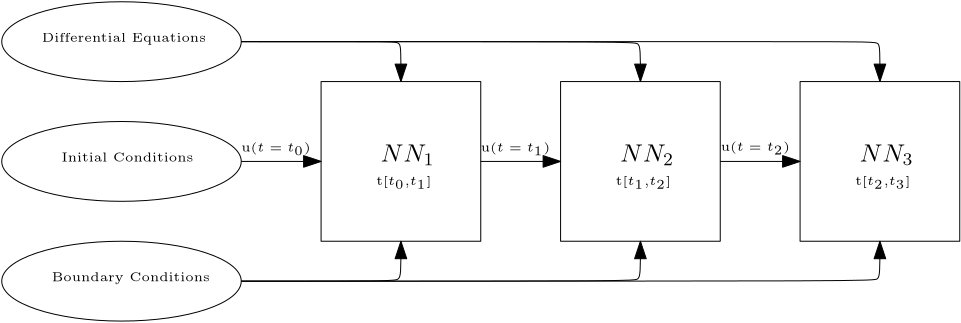}
\caption{Sequential Training Procedure}
\label{fig:seq_train}
\end{figure}
Drawing inspiration from architectures designed to handle sequential temporal data such as RNNs, LSTMs~\cite{hochreiter1997long}, the neural network is treated as a sliding window function. The network is trained for a range of time steps, and the output of the final time step of the network is taken as the initial condition for the next training phase. This autoregressive method allows for moderate function complexity while enabling a modular structure for solving large time frames. As a result, the memory requirements for training PINNs also reduce considerably for large time frames, allowing to be trained on a single GPU. The choice of time increment $\Delta t$, i.e. the length of the time frame each phase trains, is a hyper-parameter that is decided upon the number of collocation points, size of the network and the geometry of the boundary.

This form of discretization is different from time-stepping done using Euler schemes used in finite element methods, and other discretization methods for PINNs~\cite{biesek2023burgers}, as each network still contains complete time information within its range. So, to query the network at any time, one has to load the network weights corresponding to that particular time range and infer the network, thus maintaining the continuous nature of the solution function. Depending on the choice, the proposed training method can be implemented in two ways. In the first method, for each phase, a new network is initialized after the previous phase completes the training. In the second method, the one used in this paper, the same network is trained recursively in each phase. The latter method can take advantage of the principles of transfer learning\cite{9134370}, avoiding training from scratch in every phase. A schematic of the new training method is illustrated in \Cref{fig:seq_train}. The developed code can be downloaded from github repository \href{https://github.com/nexushaiku/PINN2DMovingSource}{PINN2D}.

\section{Numerical Results}
\label{sec:numexa}
\begin{figure}[hptb]
\centering 
\includegraphics[scale=0.6]{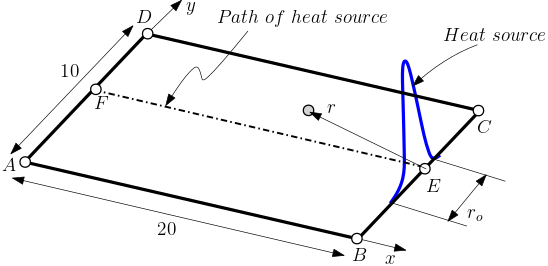}
\caption{Schematic of geometry and boundary conditions}
\label{fig:geom_des}
\end{figure}

In this section, the PINNs described in \Cref{sec:pinnoverview} are implemented in two dimensions. For the numerical analysis, the computational domain is assumed to be rectangular domain with 20 mm$\times$10 mm. The heat source is assumed to be a Gaussian with its peak at the mid-side of the shortest edge at $t=$ 0 s. The heat source is assumed to travel in a straight line $E-F$. See \Cref{fig:geom_des} for a schematic representation of the boundary and the path of the heat source. The following boundary conditions are enforced, $\forall t$ on the exterior boundaries: 
\begin{equation}
\left\{
    \begin{aligned}
        {\rm Along~A-D},~\Gamma_{AD}: &~T = 298 K \\
        {\rm Along~A-B},~\Gamma_{AB}: &~\mathbf{q}\cdot \mathbf{n} = 0.001 W \\
        {\rm Along~B-C},~\Gamma_{AB}: &~\mathbf{q}\cdot \mathbf{n}  = 0.001 W \\
        {\rm Along~C-D},~\Gamma_{CD}: &~\mathbf{q}\cdot \mathbf{n} = 0.001 W
    \end{aligned}
\right.
\end{equation}
\Cref{tab:parameters} lists the parameters used for the computations. The results from the PINNs are compared with the results from the traditional FEM. In case of PINNs, a 9-layer network deep feed-forward with 128 neurons at each layer is used. Adam optimizer is chosen along with a `$\tanh()$' activation function for the hidden layers. The network is trained for 20,000 epochs at every stage, and $\lambda_{ic}, \lambda_{bc}, \lambda_r$ are chosen to be 250,~250,~1000, respectively. 

For the FEM, the domain is discretized with 3-noded triangular elements. Based on progressive mesh refinement, a mesh size $h=$ 0.05 mm is found to be adequate and is used for spatial discretization. A backward Euler scheme with $\Delta t=$ 0.1s is used for temporal discretization. All numerical simulations were performed on Intel Xeon Ice lake 6326@2.9 G, 512 GB RAM equipped with Nvidia Tesla A100-80 GB graphics card.

\begin{table}[htpb]
\centering
\begin{tabular}{lr}
\hline
\textbf{Parameter}               & \textbf{Value}       \\ \hline
Length                           & 20 mm              \\ 
Width                            & 10 mm               \\
Specific heat                    & 658 J/kg/$\circ$         \\ 
Material density                 & $7.6 \cdot 10^{-6}$ kg/mm$^3$ \\
Thermal conductivity             & 0.025 W/mm/°C       \\ 
Maximum heat source rate         & 5 W/mm$^3$          \\ 
Radius of heat source            & 1 mm                \\
Velocity of heat source          & 2 mm/s              \\
\hline
\end{tabular}
\caption{Parameter values for the numerical simulation}
\label{tab:parameters}
\end{table}

\begin{figure}[htpb]
\subfloat[PINN $t=$ 2s]{\includegraphics[scale=0.35]{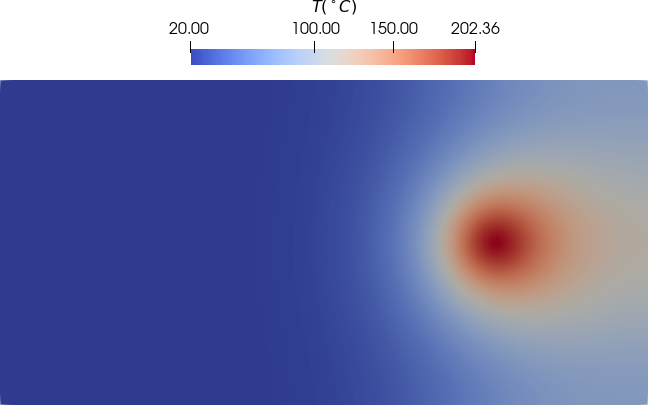}}\hspace{2cm}
\subfloat[FEM $t=$ 2s]{\includegraphics[scale=0.35]{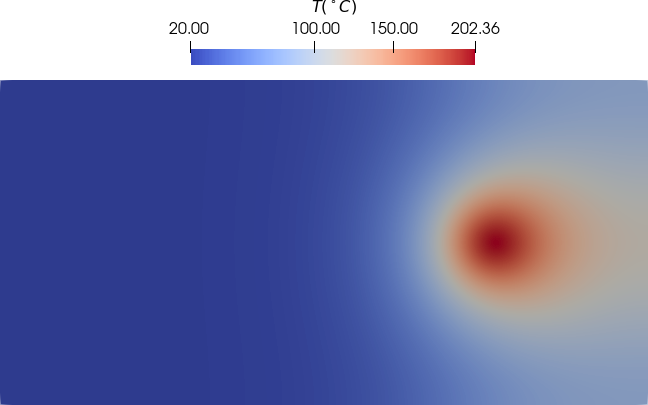}}\\
\subfloat[PINN $t=$ 8s]{\includegraphics[scale=0.35]{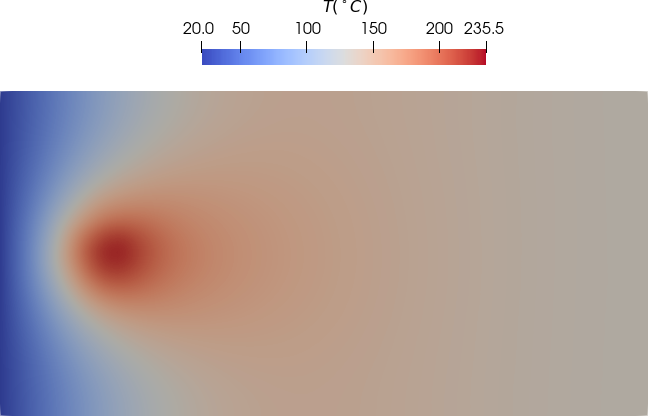}}\hspace{2cm}
\subfloat[FEM $t=$ 8s]{\includegraphics[scale=0.35]{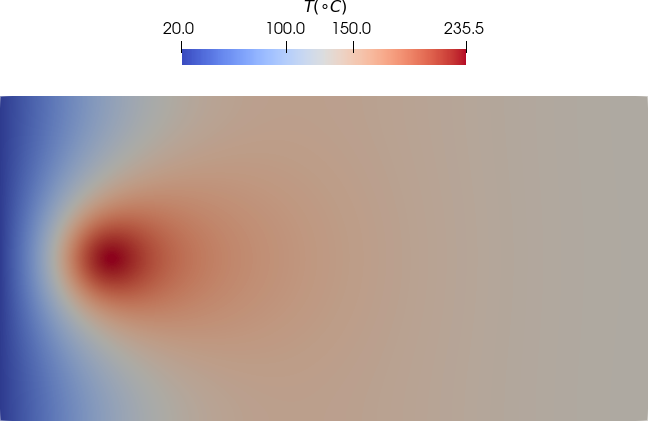}}
\caption{Comparison of temperature evolution along the path $E-F$ (c.f. \Cref{fig:geom_des}) between the traditional FEM and the PINNs}
\label{fig:contourplot_diff_time}
\end{figure}
\begin{figure}[htpb]
\subfloat[$t=$ 2s]{\includegraphics[scale=0.42]{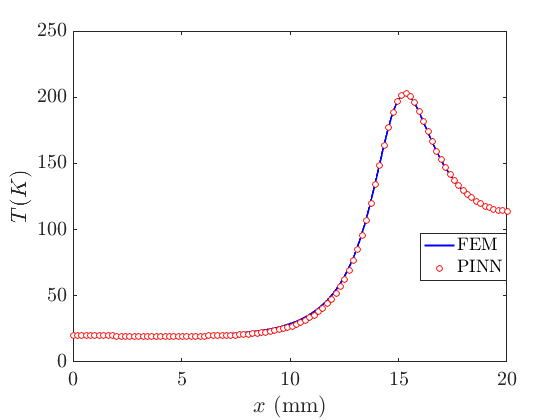}}
\subfloat[$t=$ 4s]{\includegraphics[scale=0.42]{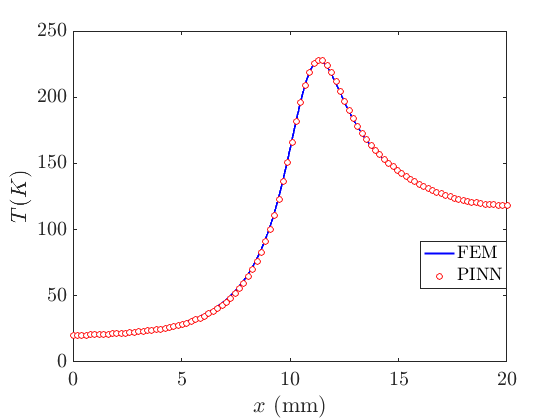}}\\
\subfloat[$t=$ 6s]{\includegraphics[scale=0.42]{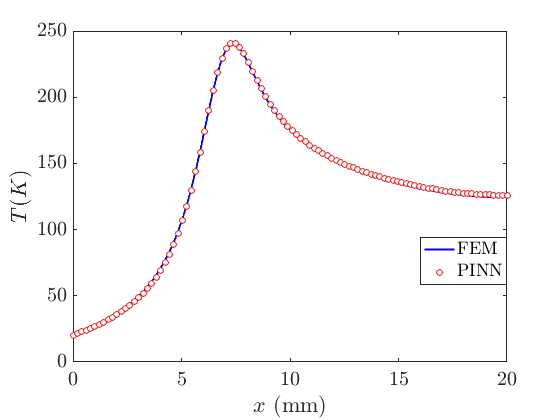}}
\subfloat[$t=$ 8s]{\includegraphics[scale=0.42]{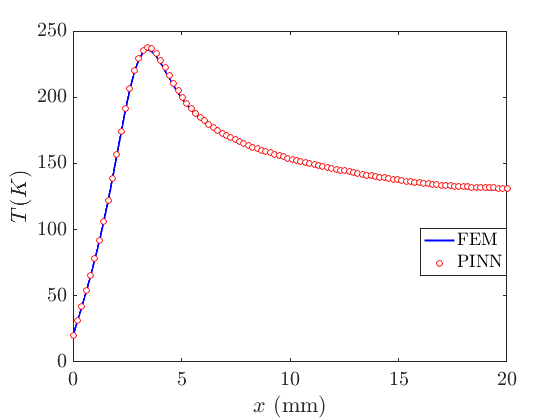}}
\caption{Comparison of temperature evolution along the path $E-F$ (c.f. \Cref{fig:geom_des}) between the traditional FEM and the PINNs}
\label{fig:lineplotcompa}
\end{figure}

\Cref{fig:contourplot_diff_time} shows the temperature distribution at $t=$ 2 s and $t=$ 8 s. The results from the present framework are compared with the FEM and a good agreement is observed. To investigate the PINNs ability to map the peak temperature at every time step, temperature evolution along the path $E-F$ (c.f. \Cref{fig:geom_des}) are computed. These values are obtained by querying the network for the PINN, and using a linear interpolator to compute values exactly, similar to that done in the FEM. \Cref{fig:lineplotcompa} shows the temperature evolution along the path $E-F$ at two different time instances, namely, $t=$ 2s and 6s. The results from the proposed framework show excellent agreement with those of the FEM.

Next, the influence of the velocity of the heat source on the temperature evolution is studied. Three different velocities are considered, that is, $v = $ (0.5, 1, 2) mm/s, and the corresponding temperature distribution along the path $E-F$ at two different time instances, that is, $t=$ 2 and 4 s, is shown in \Cref{fig:lineplotcompa_velocity}. It is seen that as the velocity of the heat source is increased for the same spatial location, the maximum temperature reached reduces, as expected. This could be attributed to the fact that with increasing velocity of the heat source, the time the region is exposed to the heat source reduces. Another aspect worth mentioning here is the training times. \revans{The PINNs with all the training phases takes $\approx$ 3100 s to complete simulation. The FEM on the other hand depends highly on mesh sizing and time step. At the fine mesh size ($h=$ 0.05 mm)and a time step of 0.1s, the FEM took 6900 s to complete the simulation.}

It is clear that PINNs demonstrate the ability to solve a moving heat source problem with mixed boundary conditions. It is worth mentioning here that apart from initial and boundary conditions, no additional data have been provided for extrapolation. The results are computed directly from the governing differential equations and corresponding initial and boundary conditions (c.f. \Cref{eq:def}). However, this also means that the network is tasked with mapping a more complex function in 3 dimensions $(x,y,t)$, and hence scaling proves to be an issue, especially since the factor that determines complexity in PINNs is geometry, instead of mesh size, as is the case with FEM. The temporal complexity manifests itself through some artifacting at $t=$ 2s which is not present in any other time frames.

\section{Conclusions}
\label{sec:concl}
In this paper, the effectiveness of PINNs in the application of moving heat source is studied. A mixed Dirichlet-Neumann boundary problem is imposed with a moving Gaussian heat source. Taking into account the complexity of the solution function, a new training regime is proposed that uses a continuous time-stepping through transfer learning. The overall time interval of the simulation is broken down into smaller intervals of $\Delta t$ period each. A single network is initialized and each phase is trained on the same network, with the initial condition of the new network the output of the final time step of the previous phase. This allows for training for extended time periods, without scaling network parameters and hyperparameters accordingly, without compromising on the resolution of the final result. Training times also scale linearly according to the number of training stages, which is a result of $\Delta t$ and overall time interval. The results of the present formulation are compared with the traditional finite element that employs 3-noded triangular elements to discretize the computational domain. From the numerical results, it is inferred that the PINNs could be seen as an alternate effective framework. To achieve similar accuracy, the PINNs require $\approx$ half the time that of the FEM. The influence of velocity of the heat source on the temperature distribution is also studied and as expected, the maximum temperature reached for the same spatial location reduces with increasing velocity. The velocity of the heat source is an important parameter that decides the temperature evolution, development of residual stresses and in turn the structural integrity. Combining the proposed framework with phase field, extension to higher-dimensions and coupled problems with temperature dependent properties will be a scope of future work.

\section*{Appendix}
The framework is implemented  can be downloaded from
\href{https://github.com/nexushaiku/PINN2DMovingSource}{https://github.com/nexushaiku/PINN2DMovingSource}

\printbibliography

\end{document}